\documentclass[11pt]{article}%
\usepackage{amsmath}
\usepackage{amsfonts}
\usepackage{amssymb}
\usepackage{graphicx}%
\providecommand{\U}[1]{\protect\rule{.1in}{.1in}}
\newtheorem{theorem}{Theorem}[section]

\numberwithin{equation}{section}
\newenvironment{Proof}[1][Proof]{\noindent\textbf{#1.} }{\ \rule{0.5em}{0.5em}}
\begin{document}

\title{Spectrum completion and inverse Sturm-Liouville problems }
\author{Vladislav V. Kravchenko\\{\small Departamento de Matem\'{a}ticas, Cinvestav, Unidad Quer\'{e}taro, }\\{\small Libramiento Norponiente \#2000, Fracc. Real de Juriquilla,
Quer\'{e}taro, Qro., 76230 MEXICO.}\\{\small e-mail: vkravchenko@math.cinvestav.edu.mx}}
\maketitle

\begin{abstract}
Given a finite set of eigenvalues of a regular Sturm-Liouville problem for the
equation $-y^{\prime\prime}+q(x)y=\lambda y$, the potential $q(x)$ of which is
unknown. We show the possibility to compute more eigenvalues without any
additional information on the potential $q(x)$. Moreover, considering the
Sturm-Liouville problem with the boundary conditions $y^{\prime}(0)-hy(0)=0$
and $y^{\prime}(\pi)+Hy(\pi)=0$, where $h$, $H$ are some constants, we
complete its spectrum without additional information neither on the potential
$q(x)$ nor on the constants $h$ and $H$. The eigenvalues are computed with a
uniform absolute accuracy. Based on this result we propose a new method for
numerical solution of the inverse Sturm-Liouville problem of recovering the
potential from two spectra. The method includes the completion of the spectra
in the first step and reduction to a system of linear algebraic equations in
the second. The potential $q(x)$ is recovered from the first component of the
solution vector. The approach is based on special Neumann series of Bessel
functions representations for solutions of Sturm-Liouville equations
possessing remarkable properties and leads to an efficient numerical algorithm
for solving inverse Sturm-Liouville problems.

\end{abstract}

\section{Introduction}

Let $q\in L_{2}(0,\pi)$ be real valued. Consider the Sturm-Liouville equation
\begin{equation}
-y^{\prime\prime}+q(x)y=\lambda y,\quad x\in(0,\pi), \label{SL equation}%
\end{equation}
where $\lambda\in\mathbb{C}$. In this work we explore the following surprising
possibility. Suppose, there are given several first eigenvalues of a
corresponding regular Sturm-Liouville problem, and the potential $q(x)$ is
unknown, we compute an arbitrarily large number of subsequent eigenvalues with
a uniform absolute accuracy.

We show that a very limited number of the eigenvalues may be sufficient for
computing hundreeds of subsequent eigenvalues with a remarkable accuracy. Of
course, it cannot be a question of calculating the eigenvalues based on their
asymptotics, because such a reduced number of known eigenvalues is definitely
not enough for obtaining the asymptotics, nor is it possible to talk about
recovering the potential from several eigenvalues of one spectrum.

Our approach is based on completely different ideas. In \cite{KNT} special
representations for solutions of (\ref{SL equation}) and for their derivatives
in the form of so-called Neumann series of Bessel functions (NSBF) were
obtained, possessing certain important properties. First of all, the
remainders of the series admit estimates independent of $\rho=\sqrt{\lambda}$
for all real $\rho$ or $\rho$ belonging to a strip $\left\vert
\operatorname{Im}\rho\right\vert <C$, where $C>0$. Simply put, the truncated
series approximate equally well the exact solutions and their first
derivatives independently of the largeness of $\left\vert \operatorname{Re}%
\rho\right\vert $. This is extremely useful when considering direct and
inverse spectral problems because it allows one to operate on large intervals
of $\rho$ and consequently of $\lambda$. Second, the knowledge of the very
first coefficient of the series is sufficient for recovering the potential
$q(x)$. These unique features of the NSBF representations were used in
\cite{KNT}, \cite{KT2018Calcolo}, \cite{KrTorbaCastillo2018} for solving
direct Sturm-Liouville problems and in \cite{Kr2019JIIP}, \cite{KrBook2020},
\cite{KKK2020}, \cite{KST2020IP}, \cite{KT2021 IP1}, \cite{KT2021 IP2} for
solving inverse Sturm-Liouville problems. In the present work we show that the
NSBF representations allow one to complete the spectrum of the Sturm-Liouville
problem and how this spectrum completion is used for solving the two-spectra
inverse Sturm-Liouville problem.

We develop the spectrum completion technique for the Dirichlet-Dirichlet
spectrum, that is the eigenvalues $\left\{  \nu_{n}\right\}  _{n=1}^{\infty}$
of (\ref{SL equation}) subject to the boundary conditions%
\begin{equation}
y(0)=y(\pi)=0, \label{DD cond}%
\end{equation}
for the Dirichlet-Neumann spectrum, that is the eigenvalues $\left\{
\lambda_{n}\right\}  _{n=0}^{\infty}$ of (\ref{SL equation}) subject to the
boundary conditions%
\begin{equation}
y(0)=y^{\prime}(\pi)=0, \label{DN cond}%
\end{equation}
as well as for the Sturm-Liouville problem with the boundary conditions
\[
y^{\prime}(0)-hy(0)=y^{\prime}(\pi)+Hy(\pi)=0,
\]
where $h$ and $H$ are unknown real constants. We complete the spectrum without
knowing the values of these constants.

Moreover, often, especially when solving inverse Sturm-Liouville problems
numerically, the authors are forced to assume that besides the boundary
conditions the important quantity
\[
\omega:=\frac{1}{2}\int_{0}^{\pi}q(t)dt
\]
is known. It appears in the asymptotics of the eigenvalues (see, e.g.,
(\ref{asymptotics DD}) and (\ref{asymptotics DN})) and is used in different
steps of most existing algorithms. To extract this parameter from the
asymptotics of the eigenvalues a considerable number of the eigenvalues are
required. Our approach allows us to compute the parameter $\omega$ from very
few eigenvalues (see Subsection \ref{Subsect DN completion} below) due to its
close relation to the first coefficient of the NSBF representation for the
derivative of the solution.

As an application of the spectrum completion technique we develop a new method
for solving the classical inverse Sturm-Liouville problem consisting in
numerical recovering the potential $q(x)$ from a finite set of eigenvalues
from two spectra.

Several methods have been proposed for numerical solution of inverse
Sturm-Liouville problems (see \cite{Brown et al 2003}, \cite{Drignei 2015},
\cite{Gao et al 2013}, \cite{Gao et al 2014}, \cite{IgnatievYurko},
\cite{Kammanee Bockman 2009}, \cite{Kr2019JIIP}, \cite{Lowe et al 1992},
\cite{Neamaty et al 2017}, \cite{Neamaty et al 2019}, \cite{Rohrl},
\cite{Rundell Sacks}, \cite{Sacks}). However, usually they require the
knowledge of additional parameters like, e.g., the parameter $\omega$.

The method proposed in the present work reduces the problem to a system of
linear algebraic equations for finding the first coefficient of the NSBF
representation for the solution. The accuracy of the method relies on the
spectrum completion technique. To the difference from the numerical methods
developed earlier on the base of the NSBF representations (see
\cite{Kr2019JIIP}, \cite{KrBook2020}, \cite{KKK2020}, \cite{KST2020IP},
\cite{KT2021 IP1}, \cite{KT2021 IP2}) here the system of linear algebraic
equations is obtained without using the Gelfand-Levitan integral equation but
a relation between the eigenfunctions normalized at the opposite endpoints.
The method is simple in its numerical realization, accurate and fast.

\section{Preliminaries}

Let us recall the asymptotics of the Dirichlet-Dirichlet and Dirichlet-Neumann
eigenvalues of equation (\ref{SL equation}). The square roots of the
eigenvalues of the Sturm-Liouville problem (\ref{SL equation}), (\ref{DD cond}%
) satisfy the asymptotic relation (see, e.g., \cite[p. 18]{Yurko2007})%
\begin{equation}
\mu_{n}=\sqrt{\nu_{n}}=n+\frac{\omega}{\pi n}+\frac{\kappa_{n}}{n}
\label{asymptotics DD}%
\end{equation}
where $\left\{  \kappa_{n}\right\}  \in\ell_{2}$ and
\begin{equation}
\omega=\frac{1}{2}\int_{0}^{\pi}q(t)dt. \label{omega}%
\end{equation}
The square roots of the eigenvalues of the Sturm-Liouville problem
(\ref{SL equation}), (\ref{DN cond}) satisfy the asymptotic relation (see,
e.g., \cite[p. 18]{Yurko2007})%
\begin{equation}
\rho_{n}=\sqrt{\lambda_{n}}=n+\frac{1}{2}+\frac{\omega}{\pi n}+\frac
{\varkappa_{n}}{n},\quad\left\{  \varkappa_{n}\right\}  \in\ell_{2}.
\label{asymptotics DN}%
\end{equation}

By $\varphi(\rho,x)$ and $S(\rho,x)$ we denote the solutions of the equation
\begin{equation}
-y^{\prime\prime}(x)+q(x)y(x)=\rho^{2}y(x),\quad x\in(0,\pi) \label{Schr}%
\end{equation}
satisfying the initial conditions
\[
\varphi(\rho,0)=1,\quad\varphi^{\prime}(\rho,0)=h,
\]%
\begin{equation}
S(\rho,0)=0,\quad S^{\prime}(\rho,0)=1, \label{init S}%
\end{equation}
where $h$ is some (complex) constant. Here $\rho=\sqrt{\lambda}$,
$\operatorname{Im}\rho\geq0$. The main tool used in the present work is the
series representations obtained in \cite{KNT} for the solutions of
(\ref{Schr}) and their derivatives.

\begin{theorem}
[\cite{KNT}]\label{Th NSBF} The solutions $\varphi(\rho,x)$ and $S(\rho,x)$
and their derivatives with respect to $x$ admit the following series
representations
\begin{align}
\varphi(\rho,x)  &  =\cos\left(  \rho x\right)  +\sum_{n=0}^{\infty}%
(-1)^{n}g_{n}(x)\mathbf{j}_{2n}(\rho x),\label{phiNSBF}\\
S(\rho,x)  &  =\frac{\sin\left(  \rho x\right)  }{\rho}+\frac{1}{\rho}%
\sum_{n=0}^{\infty}(-1)^{n}s_{n}(x)\mathbf{j}_{2n+1}(\rho x),\label{S}\\
\varphi^{\prime}(\rho,x)  &  =-\rho\sin\left(  \rho x\right)  +\left(
h+\frac{1}{2}\int_{0}^{x}q(t)\,dt\right)  \cos\left(  \rho x\right)
+\sum_{n=0}^{\infty}(-1)^{n}\gamma_{n}(x)\mathbf{j}_{2n}(\rho
x),\label{phiprimeNSBF}\\
S^{\prime}(\rho,x)  &  =\cos\left(  \rho x\right)  +\frac{1}{2\rho}\left(
\int_{0}^{x}q(t)\,dt\right)  \sin\left(  \rho x\right)  +\frac{1}{\rho}%
\sum_{n=0}^{\infty}(-1)^{n}\sigma_{n}(x)\mathbf{j}_{2n+1}(\rho x),
\label{Sprime}%
\end{align}
where $\mathbf{j}_{k}(z)$ stands for the spherical Bessel function of order
$k$ (see, e.g., \cite{AbramowitzStegunSpF}). The coefficients $g_{n}(x)$,
$s_{n}(x)$, $\gamma_{n}(x)$ and $\sigma_{n}(x)$ can be calculated following a
simple recurrent integration procedure (see \cite{KNT} or \cite[Sect.
9.4]{KrBook2020}), starting with
\begin{align}
g_{0}(x)  &  =\varphi(0,x)-1,\quad s_{0}(x)=3\left(  \frac{S(0,x)}%
{x}-1\right)  ,\label{beta0}\\
\gamma_{0}(x)  &  =g_{0}^{\prime}(x)-h-\frac{1}{2}\int_{0}^{x}q(t)\,dt,\quad
\sigma_{0}(x)=\frac{s_{0}(x)}{x}+s_{0}^{\prime}(x)-\frac{3}{2}\int_{0}%
^{x}q(t)\,dt.\nonumber
\end{align}
For every $\rho\in\mathbb{C}$ all the series converge pointwise. For every
$x\in\left[  0,\pi\right]  $ the series converge uniformly on any compact set
of the complex plane of the variable $\rho$, and the remainders of their
partial sums admit estimates independent of $\operatorname{Re}\rho$.
\end{theorem}

This last feature of the series representations (the independence of
$\operatorname{Re}\rho$ of the estimates for the remainders) is a direct
consequence of the fact that the representations are obtained by expanding the
integral kernels of the transmutation operators (for their theory we refer to
\cite{LevitanInverse}, \cite{Marchenko}, \cite{SitnikShishkina Elsevier}) into
Fourier-Legendre series (see \cite{KNT} and \cite[Sect. 9.4]{KrBook2020}). It
is of crucial importance for what follows. In particular, it means that for
$S_{N}(\rho,x):=\frac{\sin\left(  \rho x\right)  }{\rho}+\frac{1}{\rho}%
\sum_{n=0}^{N}(-1)^{n}s_{n}(x)\mathbf{j}_{2n+1}(\rho x)$ (and analogously for
$S_{N}^{\prime}(\rho,x)$) the estimate holds
\begin{equation}
\left\vert S(\rho,x)-S_{N}(\rho,x)\right\vert <\varepsilon_{N}(x)
\label{estim S}%
\end{equation}
for all $\rho\in\mathbb{R}$, where $\varepsilon_{N}(x)$ is a positive function
tending to zero when $N\rightarrow\infty$. That is, the approximate solution
$S_{N}(\rho,x)$ approximates the exact one equally well for small and for
large values of $\rho$. This is especially convenient when considering direct
and inverse spectral problems. Moreover, for a fixed $z$ the numbers
$\mathbf{j}_{k}(z)$ rapidly decrease as $k\rightarrow\infty$, see, e.g.,
\cite[(9.1.62)]{AbramowitzStegunSpF}. Hence, the convergence rate of the
series for any fixed $\rho$ is, in fact, exponential.

More detailed estimates for the series remainders depending on the regularity
of the potential can be found in \cite{KNT}.

Note that formulas (\ref{beta0}) indicate that the potential $q(x)$ can be
recovered from the first coefficients of the series (\ref{phiNSBF}) or
(\ref{S}). We have
\begin{equation}
q(x)=\frac{g_{0}^{\prime\prime}(x)}{g_{0}(x)+1} \label{qi from g0}%
\end{equation}
and
\begin{equation}
q(x)=\frac{\left(  xs_{0}(x)\right)  ^{\prime\prime}}{xs_{0}(x)+3x}.
\label{qi from s0}%
\end{equation}

Note that the square roots of the Dirichlet-Dirichlet eigenvalues coincide
with zeros of the function $S(\rho,\pi)$:%
\[
S(\mu_{n},\pi)=0,\quad n=1,2,\ldots,
\]
while the square roots of the Dirichlet-Neumann eigenvalues coincide with
zeros of the function $S^{\prime}(\rho,\pi)$:%
\[
S^{\prime}(\rho_{n},\pi)=0,\quad n=0,1,\ldots.
\]

\section{Spectrum completion}

\subsection{Dirichlet-Dirichlet spectrum}

Given several first Dirichlet-Dirichlet eigenvalues $\left\{  \mu_{n}%
^{2}\right\}  _{n=1}^{N_{1}}$, let us use them to calculate the first
coefficients $s_{0}(\pi)$, $s_{1}(\pi)$,...,$s_{N}(\pi)$, where $N\leq
N_{1}+1$. First, it is convenient to consider the shifted potential
\[
\widetilde{q}(x):=q(x)-\mu_{1}^{2}.
\]
That is, instead of the problem (\ref{SL equation}), (\ref{DD cond}) we
consider the problem
\begin{equation}
-y^{\prime\prime}(x)+\widetilde{q}(x)y(x)=\rho^{2}y(x),\quad x\in(0,\pi),
\label{Schrqtilde}%
\end{equation}%
\begin{equation}
y(0)=y(\pi)=0. \label{DirProblem2}%
\end{equation}
Obviously, its eigenfunctions do not change while the eigenvalues are shifted:%
\[
\widetilde{\mu}_{1}^{2}=0,\quad\widetilde{\mu}_{2}^{2}=\mu_{2}^{2}-\mu_{1}%
^{2},\quad\widetilde{\mu}_{3}^{2}=\mu_{3}^{2}-\mu_{1}^{2},\ldots.
\]
The solution of (\ref{Schrqtilde}) satisfying the initial conditions
(\ref{init S}) we denote as $\widetilde{S}(\rho,x)$. Its NSBF representation
has the form
\begin{equation}
\widetilde{S}(\rho,x)=\frac{\sin\left(  \rho x\right)  }{\rho}+\frac{1}{\rho
}\sum_{n=0}^{\infty}(-1)^{n}\widetilde{s}_{n}(x)\mathbf{j}_{2n+1}(\rho x),
\label{Stilde}%
\end{equation}
where the coefficients $\widetilde{s}_{n}(x)$, in general, do not coincide
with the coefficients $s_{n}(x)$, however it is clear that
\[
S(\rho,x)=\widetilde{S}(\sqrt{\rho^{2}+\mu_{1}^{2}},x).
\]
Note that
\begin{equation}
\widetilde{S}(0,\pi)=0, \label{Stilde0pi=0}%
\end{equation}
since zero is a Dirichlet-Dirichlet eigenvalue of $\widetilde{q}(x)$.

On the other hand, from (\ref{Stilde}) we have
\begin{equation}
\widetilde{S}(0,x)=x+\frac{x\widetilde{s}_{0}(x)}{3}, \label{Stilde0}%
\end{equation}
where we take into account that $\mathbf{j}_{1}(z)\sim\frac{z}{3}$,
$z\rightarrow0$ and more generally,
\[
\mathbf{j}_{n}(z)\sim\frac{z^{n}}{\left(  2n+1\right)  !!},\quad
z\rightarrow0.
\]
Substituting $x=\pi$ into (\ref{Stilde0}) and taking into account
(\ref{Stilde0pi=0}), we obtain%
\begin{equation}
\widetilde{s}_{0}(\pi)=-3. \label{s0tildepi}%
\end{equation}

Several subsequent coefficients $\widetilde{s}_{n}(\pi)$, $n=1,\ldots,N$ can
be found from the equations%
\begin{equation}
\widetilde{S}_{N}(\widetilde{\mu}_{k},\pi)=0,\quad k=2,3,\ldots,N_{1},
\label{StildeN=0}%
\end{equation}
where%
\begin{equation}
\widetilde{S}_{N}(\rho,x)=\frac{\sin\left(  \rho x\right)  }{\rho}+\frac
{1}{\rho}\sum_{n=0}^{N}(-1)^{n}\widetilde{s}_{n}(x)\mathbf{j}_{2n+1}(\rho x).
\label{StildeN}%
\end{equation}
From (\ref{StildeN=0}) we obtain the system of linear algebraic equations for
the coefficients $\widetilde{s}_{n}(\pi)$:%
\[
\sum_{n=0}^{N}(-1)^{n}\widetilde{s}_{n}(\pi)\mathbf{j}_{2n+1}(\widetilde{\mu
}_{k}\pi)=-\sin\left(  \widetilde{\mu}_{k}\pi\right)  ,\quad k=2,3,\ldots
,N_{1}.
\]
Taking into account (\ref{s0tildepi}), we obtain
\begin{equation}
\sum_{n=1}^{N}(-1)^{n}\widetilde{s}_{n}(\pi)\mathbf{j}_{2n+1}(\widetilde{\mu
}_{k}\pi)=3\mathbf{j}_{1}(\widetilde{\mu}_{k}\pi)-\sin\left(  \widetilde{\mu
}_{k}\pi\right)  ,\quad k=2,3,\ldots,N_{1}. \label{system sn tilde}%
\end{equation}
Solving this system of equations we compute $\widetilde{s}_{1}(\pi
),\ldots,\widetilde{s}_{N}(\pi)$. This gives us the possibility to compute an
arbitrary number of the Dirichlet-Dirichlet eigenvalues of the problem
(\ref{Schrqtilde}), (\ref{DirProblem2}) and consequently of the original
problem (\ref{SL equation}), (\ref{DD cond}). Indeed, the function%
\begin{equation}
\widetilde{S}_{N}(\rho,\pi)=\frac{\sin\left(  \rho\pi\right)  }{\rho}+\frac
{1}{\rho}\sum_{n=0}^{N}(-1)^{n}\widetilde{s}_{n}(\pi)\mathbf{j}_{2n+1}(\rho
\pi) \label{StildeNatpi}%
\end{equation}
approximates the solution $\widetilde{S}(\rho,\pi)$ at $x=\pi$ for any value
of $\rho$. Moreover, for $\rho\in\mathbb{R}$ we have the estimate
\[
\left\vert \widetilde{S}(\rho,\pi)-\widetilde{S}_{N}(\rho,\pi)\right\vert
<\widetilde{\varepsilon}_{N},
\]
where $\widetilde{\varepsilon}_{N}$ is independent of $\rho$. With the aid of
complex analysis tools the following theorem is proved.

\begin{theorem}
\label{Th closeness of zeros}For any $\varepsilon>0$ there exists such
$N\in\mathbb{N}$ that all zeros of the function $\widetilde{S}(\rho,\pi)$ are
approximated by corresponding zeros of the function $\widetilde{S}_{N}%
(\rho,\pi)$ with errors uniformly bounded by $\varepsilon$, and $\widetilde
{S}_{N}(\rho,\pi)$ has no other zeros.
\end{theorem}

\begin{Proof}
The proof of this statement is completely analogous to the proof of
Proposition 7.1 in \cite{KT2015JCAM} and consist in the use of properties of
characteristic functions of regular Sturm-Liouville problems and application
of the Rouch\'{e} theorem.
\end{Proof}

Thus, zeros of the function $\widetilde{S}_{N}(\rho,\pi)$ give us approximate
numbers $\widetilde{\mu}_{k}$ for $k=N_{1}+1,\ldots$.

\textbf{Example 1. }Consider the potential $q(x)=e^{x}$ (first Paine's test,
see \cite{PryceBook}). Then $\mu_{1}^{2}\approx4.89666937996$ (here and below
we used the Matslise package \cite{Ledoux et al} to compute the
\textquotedblleft exact\textquotedblright\ eigenvalues). In Fig. 1 the
absolute and relative errors of $\mu_{k}$ computed for $k=6,\ldots,300$ are
presented. Here five Dirichlet-Dirichlet eigenvalues were given ($N_{1}=5$),
and we present the \textquotedblleft completed\textquotedblright\ part of the
Dirichlet-Dirichlet spectrum computed with $N=4$ (so that four coefficients
$\widetilde{s}_{n}(\pi)$ in (\ref{StildeNatpi}) are computed from
(\ref{system sn tilde}) and together with the coefficient (\ref{s0tildepi})
they are used to compute the eigenvalues by finding zeros of the function
$\widetilde{S}_{N}(\rho,\pi)$). It is worth mentioning that when dealing with
the truncated systems of linear algebraic equations we do not seek to work
with the square systems. In computations a least-squares solution of an
overdetermined system (provided by Matlab, which we used in this work) gives
very satisfactory results.%

\begin{figure}
[ptb]
\begin{center}
\includegraphics[
height=3.8563in,
width=5.1304in
]%
{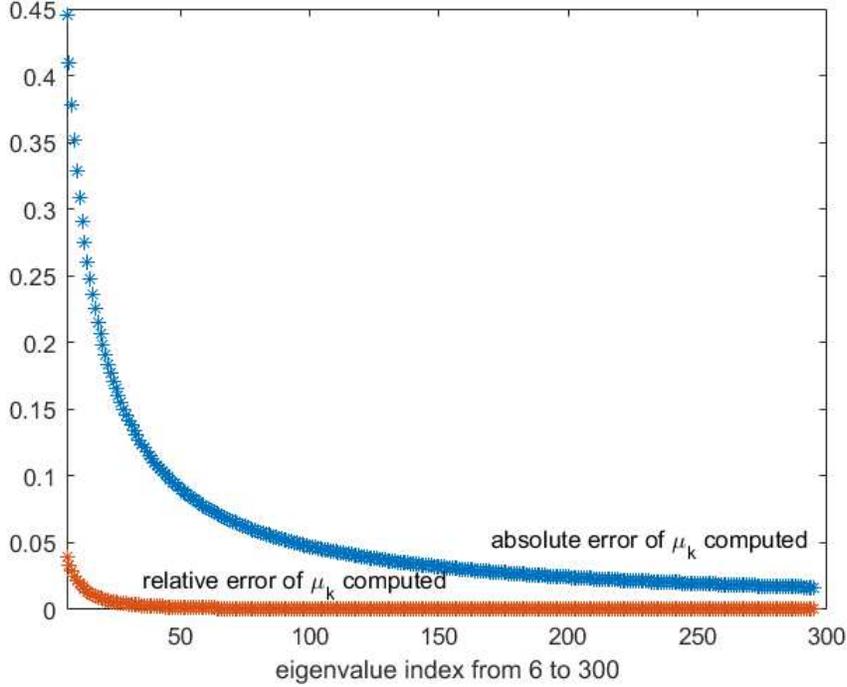}%
\caption{Absolute and relative errors of $\mu_{k}$, $k=6,\ldots,300$ computed
from five eigenvalues given for the potential $q(x)=e^{x}$.}%
\label{Fig1}%
\end{center}
\end{figure}

Of course, five eigenvalues are not enough to obtain from
(\ref{asymptotics DD}) a meaningful value of the parameter $\omega$, so it is
not clear how one could complete the spectrum in another way. Even the
knowledge of ten eigenvalues still do not give the possibility to find
$\omega$ with a reasonable accuracy. Fig. 2 shows the numbers
\begin{equation}
c_{k}=\pi k(\mu_{k}-k)-\omega,\quad k=3,\ldots,10. \label{ck}%
\end{equation}
which according to (\ref{asymptotics DD}) represent an $\ell_{2}$-convergent
sequence. The value of $c_{10}$ is approximately $-0.06$ which still differs
from zero considerably. The attempt to compute $\omega$ by minimizing the
$\ell_{2}$-norm of the sequence $\left\{  c_{k}\right\}  _{k=1}^{10}$ leads to
a large error.%

\begin{figure}
[ptb]
\begin{center}
\includegraphics[
height=3.1465in,
width=4.1868in
]%
{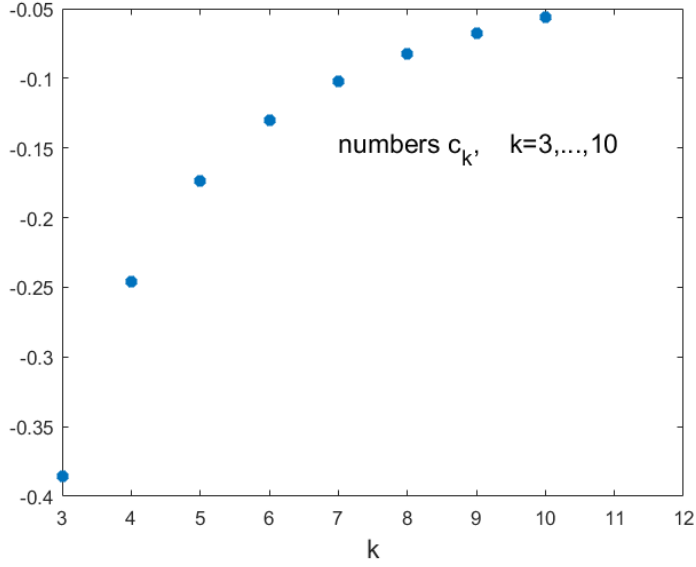}%
\caption{Numbers $c_{k}$ defined by (\ref{ck}) computed for Example 1,
$k=3,\ldots,10$.}%
\label{Fig2}%
\end{center}
\end{figure}

Moreover, even the knowledge of the exact value of the parameter $\omega$
leads to less accurate results in comparison with the spectrum completion
technique based on solving the system (\ref{system sn tilde}) and finding
zeros of the function $\widetilde{S}_{N}(\rho,\pi)$. In Fig.3 we compare the
absolute error of the numbers $\mu_{k}$ approximated by the asymptotic
relation $\mu_{k}\approx k+\frac{\omega}{\pi k}$ (where the exact value of
$\omega$ is used) with the absolute error of the values of $\mu_{k}$ obtained
with the aid of the spectrum completion technique. Even for $k=40$ the
asymptotic approximation gives a much less accurate result than the spectrum
completion technique: $4.6\cdot10^{-5}$ against $6.8\cdot10^{-7}$.

\bigskip%

\begin{figure}
[ptb]
\begin{center}
\includegraphics[
height=3.3627in,
width=4.4746in
]%
{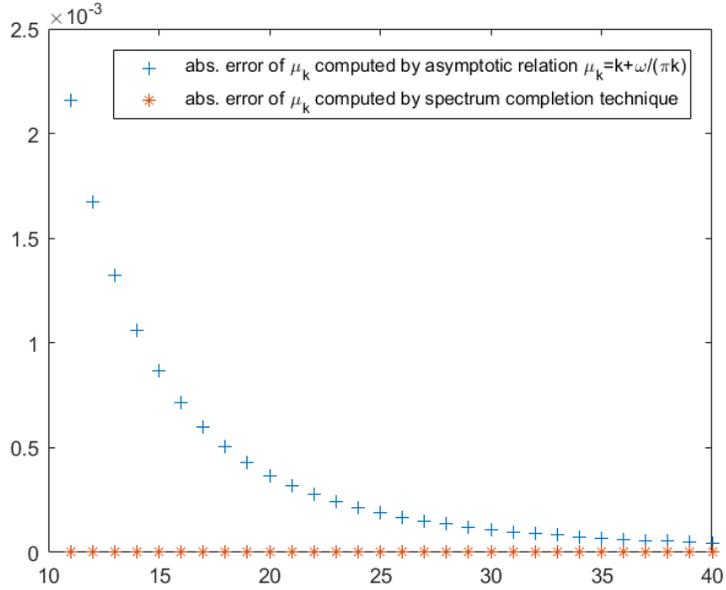}%
\caption{Even when the parameter $\omega$ is known, the spectrum completion
technique often gives more accurate results than the use of the asymptotic
relation $\mu_{k}\approx k+\frac{\omega}{\pi k}$. Here we compare the accuracy
of the \textquotedblleft asymptotic\textquotedblright\ $\mu_{k}$,
$k=11,\ldots,40$ vs. those computed by the spectrum completion technique from
ten eigenvalues given. }%
\label{Fig3}%
\end{center}
\end{figure}

A similar situation is observed in other examples.

\textbf{Example 2. }\ Consider the potential $q(x)=\frac{1}{(x+0.1)^{2}}$
(second Paine's test, see \cite{PryceBook}). In Fig. 4 we make the same
comparison: the absolute error of the numbers $\mu_{k}$ approximated by the
asymptotic relation $\mu_{k}\approx k+\frac{\omega}{\pi k}$ with the exact
value of $\omega$ being used vs. the absolute error of the values of $\mu_{k}$
obtained with the aid of the spectrum completion technique.%

\begin{figure}
[ptb]
\begin{center}
\includegraphics[
height=3.6113in,
width=4.806in
]%
{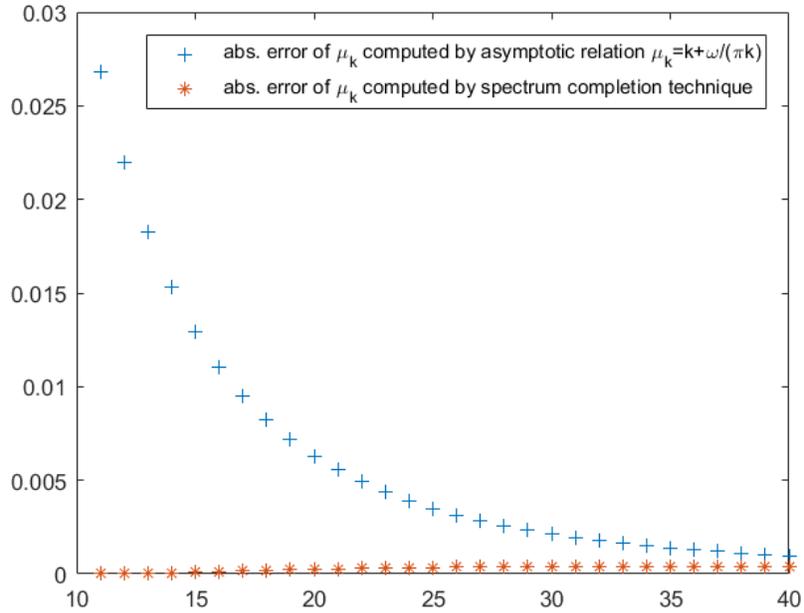}%
\caption{For the potential from Example 2 we compare the accuracy of the
\textquotedblleft asymptotic\textquotedblright\ $\mu_{k}$, $k=11,\ldots,40$
vs. those computed by the spectrum completion technique from ten eigenvalues
given. }%
\label{Fig4}%
\end{center}
\end{figure}

\subsection{Dirichlet-Neumann spectrum\label{Subsect DN completion}}

Given several first Dirichlet-Neumann eigenvalues $\left\{  \rho_{n}%
^{2}\right\}  _{n=0}^{N_{2}}$, let us use them to calculate the first
coefficients $\sigma_{0}(\pi)$, $\sigma_{1}(\pi)$,...,$\sigma_{N}(\pi)$, where
$N\leq N_{2}-1$. Similarly to the case of the Dirichlet-Dirichlet spectrum, it
is convenient to consider the shifted potential%
\[
\widehat{q}(x):=q(x)-\rho_{0}^{2}.
\]
That is, instead of the problem (\ref{SL equation}), (\ref{DN cond}) we
consider the problem%
\begin{equation}
-y^{\prime\prime}(x)+\widehat{q}(x)y(x)=\rho^{2}y(x),\quad x\in(0,\pi),
\label{Schrqhat}%
\end{equation}%
\begin{equation}
y(0)=y^{\prime}(\pi)=0. \label{DN cond 2}%
\end{equation}
Again the eigenfunctions do not change, and the eigenvalues are shifted:%
\[
\widehat{\rho}_{0}^{2}=0,\quad\widehat{\rho}_{1}^{2}=\rho_{1}^{2}-\rho_{0}%
^{2},\quad\widehat{\rho}_{2}^{2}=\rho_{2}^{2}-\rho_{0}^{2},\ldots.
\]
The solution of (\ref{Schrqhat}) satisfying the initial conditions
(\ref{init S}) we denote as $\widehat{S}(\rho,x)$. We have the relation%
\[
\widehat{S}(\rho,x)=\widetilde{S}(\sqrt{\rho^{2}+\rho_{0}^{2}-\mu_{1}^{2}%
},x).
\]
We are interested in the derivative of $\widehat{S}(\rho,x)$. It has the NSBF
representation%
\[
\widehat{S}^{\prime}(\rho,x)=\cos\left(  \rho x\right)  +\frac{1}{2\rho
}\left(  \int_{0}^{x}\widehat{q}(t)\,dt\right)  \sin\left(  \rho x\right)
+\frac{1}{\rho}\sum_{n=0}^{\infty}(-1)^{n}\widehat{\sigma}_{n}(x)\mathbf{j}%
_{2n+1}(\rho x).
\]
Note that
\[
\widehat{S}^{\prime}(0,x)=1+\frac{x}{2}\int_{0}^{x}\widehat{q}(t)\,dt+\frac
{x}{3}\widehat{\sigma}_{0}(x).
\]
Since zero is an eigenvalue of the problem (\ref{Schrqhat}), (\ref{DN cond 2}%
), we have $\widehat{S}^{\prime}(0,\pi)=0$ and thus
\[
1+\pi\widehat{\omega}+\frac{\pi}{3}\widehat{\sigma}_{0}(\pi)=0,
\]
where
\[
\widehat{\omega}:=\frac{1}{2}\int_{0}^{\pi}\widehat{q}(t)\,dt=\omega-\frac
{\pi\rho_{0}^{2}}{2}.
\]
Thus,
\begin{equation}
\widehat{\omega}=-\frac{\widehat{\sigma}_{0}(\pi)}{3}-\frac{1}{\pi}.
\label{omegahat}%
\end{equation}
Now we complete the spectrum of the problem (\ref{Schrqhat}), (\ref{DN cond 2}%
) and hence also the Dirichlet-Neumann spectrum of the potential $q(x)$. For
this we consider the equations%
\[
\widehat{S}_{N}^{\prime}(\widehat{\rho}_{k},\pi)=0,\quad k=1,2,\ldots,N_{2},
\]
which can be written in the form%
\[
\widehat{\rho}_{k}\cos\left(  \widehat{\rho}_{k}\pi\right)  +\widehat{\omega
}\sin\left(  \widehat{\rho}_{k}\pi\right)  +\sum_{n=0}^{N}(-1)^{n}%
\widehat{\sigma}_{n}(\pi)\mathbf{j}_{2n+1}(\widehat{\rho}_{k}\pi)=0.
\]
Taking into account (\ref{omegahat}), we write these equations in the form of
the system of linear algebraic equations for the coefficients $\widehat
{\sigma}_{0}(\pi)$,...,$\widehat{\sigma}_{N}(\pi)$:%
\begin{equation}
\left(  \mathbf{j}_{1}(\widehat{\rho}_{k}\pi)-\frac{\sin\left(  \widehat{\rho
}_{k}\pi\right)  }{3}\right)  \widehat{\sigma}_{0}(\pi)+\sum_{n=1}^{N}%
(-1)^{n}\widehat{\sigma}_{n}(\pi)\mathbf{j}_{2n+1}(\widehat{\rho}_{k}%
\pi)=-\widehat{\rho}_{k}\cos\left(  \widehat{\rho}_{k}\pi\right)  +\frac
{\sin\left(  \widehat{\rho}_{k}\pi\right)  }{\pi},\quad k=1,\ldots,N_{2}.
\label{syst sigman}%
\end{equation}
Solving this system we find $\widehat{\sigma}_{0}(\pi)$,...,$\widehat{\sigma
}_{N}(\pi)$ as well as the parameter $\widehat{\omega}$ (from (\ref{omegahat}%
)) and the parameter $\omega=\widehat{\omega}+\frac{\pi\rho_{0}^{2}}{2}$.

Having computed the coefficients $\widehat{\sigma}_{0}(\pi)$,...,$\widehat
{\sigma}_{N}(\pi)$ we consider the function%
\begin{equation}
\widehat{S}_{N}^{\prime}(\rho,\pi)=\cos\left(  \rho\pi\right)  +\frac
{\widehat{\omega}\sin\left(  \rho\pi\right)  }{\rho}+\frac{1}{\rho}\sum
_{n=0}^{N}(-1)^{n}\widehat{\sigma}_{n}(\pi)\mathbf{j}_{2n+1}(\rho\pi),
\label{Sprimehat}%
\end{equation}
which approximates the derivative $\widehat{S}^{\prime}(\rho,\pi)$ in such a
way that
\[
\left\vert \widehat{S}^{\prime}(\rho,\pi)-\widehat{S}_{N}^{\prime}(\rho
,\pi)\right\vert <\widehat{\varepsilon}_{N},\quad\rho\in\mathbb{R},
\]
where $\widehat{\varepsilon}_{N}$ is a positive constant. Theorem
\ref{Th closeness of zeros} is valid as well if instead of the functions
$\widetilde{S}(\rho,\pi)$ and $\widetilde{S}_{N}(\rho,\pi)$ one considers the
functions $\widehat{S}^{\prime}(\rho,\pi)$ and $\widehat{S}_{N}^{\prime}%
(\rho,\pi)$. Thus zeros of $\widehat{S}_{N}^{\prime}(\rho,\pi)$ approximate
the Dirichlet-Neumann eigenvalues of the potential $\widehat{q}(x)$.

In Fig. 5 the absolute and relative errors of $\rho_{k}$ computed for
$k=5,\ldots,300$ are presented in the case of the potential from Example 2.
Here five Dirichlet-Neumann eigenvalues were given ($N_{2}=4$), and we present
the \textquotedblleft completed\textquotedblright\ part of the
Dirichlet-Neumann spectrum computed with $N=3$ (so that four coefficients
$\widehat{\sigma}_{n}(\pi)$ in (\ref{Sprimehat}) are computed from
(\ref{syst sigman}) as well as $\widehat{\omega}$ from (\ref{omegahat})). The
absolute error of the computed $\widehat{\omega}$ was $0.8$, which is a
satisfactory result taking into account the limited number of the eigenvalues
given and that the relative error was approximately $0.19$.
\begin{figure}
[ptb]
\begin{center}
\includegraphics[
height=3.1561in,
width=4.1991in
]%
{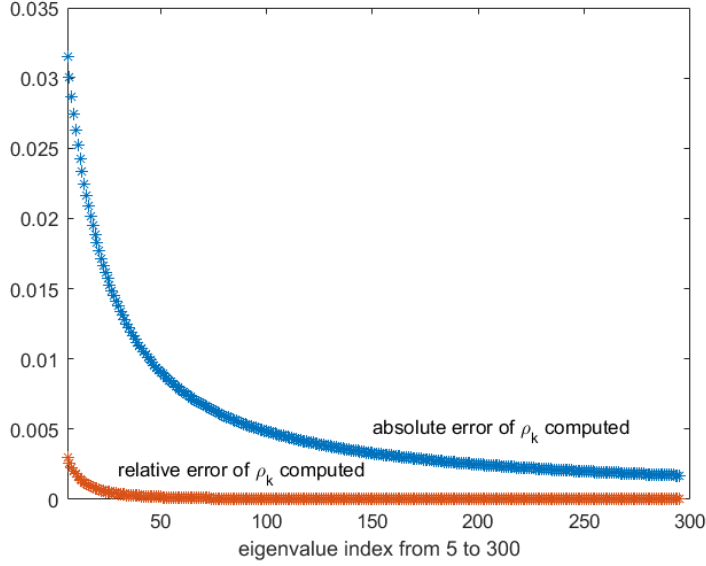}%
\caption{Absolute and relative errors of $\rho_{k}$ computed for
$k=5,\ldots,300$ are presented in the case of the potential from Example 2.
Here five Dirichlet-Neumann eigenvalues were given ($N_{2}=4$), so the
\textquotedblleft completed\textquotedblright\ part of the Dirichlet-Neumann
spectrum was computed with four coefficients $\widehat{\sigma}_{n}(\pi)$ in
(\ref{Sprimehat}) computed from (\ref{syst sigman}) and $\widehat{\omega}$
from (\ref{omegahat}).}%
\label{Fig5}%
\end{center}
\end{figure}
Naturally the accuracy improves when a larger number of the eigenvalues are
known. In Fig. 6 the results are presented for the same potential but in the
case when ten eigenvalues are known ($N_{2}=9$). Here $N=8$. The absolute
error of the computed $\widehat{\omega}$ was already $0.09$.
\begin{figure}
[ptb]
\begin{center}
\includegraphics[
height=3.2023in,
width=4.2601in
]%
{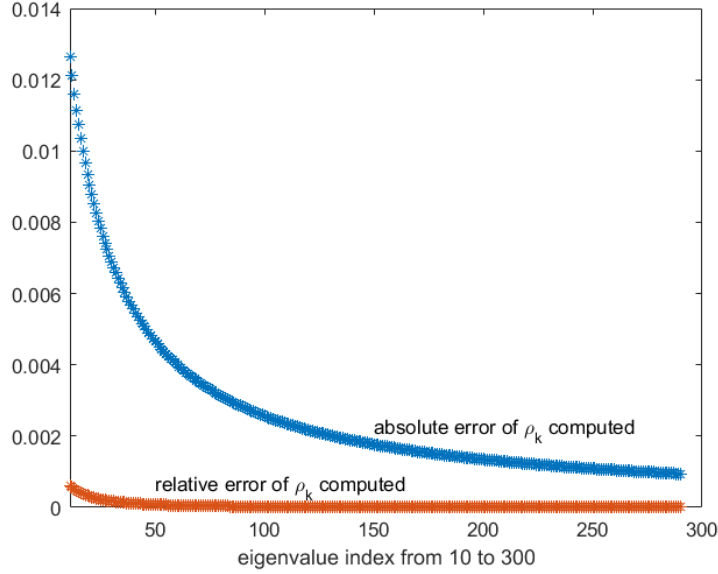}%
\caption{Same as the previous figure, but for ten Dirichlet-Neumann
eigenvalues given.}%
\label{Fig6}%
\end{center}
\end{figure}

\newpage

\subsection{Completion of other spectra}

Suppose several eigenvalues $\rho_{0}^{2}$, $\rho_{1}^{2}$,..., $\rho_{N_{3}%
}^{2}$ of the following Sturm-Liouville problem are given%
\begin{align}
-y^{\prime\prime}(x)+q(x)y(x)  &  =\rho^{2}y(x),\quad x\in(0,\pi
),\label{SLeq3}\\
y^{\prime}(0)-hy(0)  &  =0,\quad y^{\prime}(\pi)+Hy(\pi)=0, \label{SL3}%
\end{align}
where $h$ and $H$ are unknown real constants. Again, as in the previous cases,
we always can shift the eigenvalues in such a way that the first shifted
eigenvalue becomes zero. So, without loss of generality we assume that
$\rho_{0}=0$. Note that
\[
\varphi(0,x)=1+g_{0}(x)
\]
and%
\[
\varphi^{\prime}(0,x)=h+\frac{1}{2}\int_{0}^{x}q(t)dt+\gamma_{0}(x).
\]
Since zero is an eigenvalue, we have that
\[
\varphi^{\prime}(0,\pi)+H\varphi(0,\pi)=0.
\]
Thus,%
\begin{equation}
h+H+\omega=-\gamma_{0}(\pi)-Hg_{0}(\pi), \label{relation h0}%
\end{equation}
where $\omega=\frac{1}{2}\int_{0}^{\pi}q(t)dt$.

Now consider the characteristic function of the Sturm-Liouville problem
(\ref{SLeq3}), (\ref{SL3}). Taking into account (\ref{phiNSBF}) and
(\ref{phiprimeNSBF}) it can be written in the form%
\begin{align*}
\Phi(\rho)  &  :=\varphi^{\prime}(\rho,\pi)+H\varphi(\rho,\pi)\\
&  =-\rho\sin\left(  \rho\pi\right)  +\left(  h+H+\omega\right)  \cos\left(
\rho\pi\right)  +\sum_{n=0}^{\infty}(-1)^{n}\gamma_{n}(\pi)\mathbf{j}%
_{2n}(\rho\pi)+H\sum_{n=0}^{\infty}(-1)^{n}g_{n}(\pi)\mathbf{j}_{2n}(\rho\pi).
\end{align*}
Denote
\[
h_{n}:=\gamma_{n}(\pi)+Hg_{n}(\pi),\quad n=0,1,\ldots.
\]
Then, taking into account (\ref{relation h0}), we obtain%
\[
\varphi^{\prime}(\rho,\pi)+H\varphi(\rho,\pi)=-\rho\sin\left(  \rho\pi\right)
+h_{0}\left(  \mathbf{j}_{0}(\rho\pi)-\cos\left(  \rho\pi\right)  \right)
+\sum_{n=1}^{\infty}(-1)^{n}h_{n}\mathbf{j}_{2n}(\rho\pi).
\]
Several given eigenvalues $\rho_{1}^{2}$,..., $\rho_{N_{3}}^{2}$ allow us to
compute several constants $h_{n}$, $n=0,\ldots,N$, where $N\leq N_{3}-1$, from
the system of linear algebraic equations
\[
h_{0}\left(  \mathbf{j}_{0}(\rho_{k}\pi)-\cos\left(  \rho_{k}\pi\right)
\right)  +\sum_{n=1}^{N}(-1)^{n}h_{n}\mathbf{j}_{2n}(\rho_{k}\pi)=\rho_{k}%
\sin\left(  \rho_{k}\pi\right)  ,\quad k=1,\ldots,N_{3}.
\]
Thus we obtain an approximate characteristic function of the problem
(\ref{SLeq3}), (\ref{SL3})%
\[
\Phi_{N}(\rho):=h_{0}\left(  \mathbf{j}_{0}(\rho\pi)-\cos\left(  \rho
\pi\right)  \right)  +\sum_{n=1}^{N}(-1)^{n}h_{n}\mathbf{j}_{2n}(\rho\pi
)-\rho\sin\left(  \rho\pi\right)  ,
\]
whose zeros approximate the square roots of the eigenvalues of the problem
(\ref{SLeq3}), (\ref{SL3}).

We emphasize that to complete the spectrum of (\ref{SLeq3}), (\ref{SL3}) we
require no information neither on the potential $q(x)$ nor on the boundary
conditions (the constants $h$ and $H$ remain unknown). Moreover, the parameter
$\overline{\omega}:=h+H+\omega$ which appears in the second term of the
asymptotics for $\rho_{k}$ in this problem%
\[
\rho_{k}=k+\frac{\overline{\omega}}{\pi k}+\frac{\varkappa_{k}}{k}%
,\quad\left\{  \varkappa_{k}\right\}  \in\ell_{2}%
\]
is computed as well, since due to (\ref{relation h0}), $\overline{\omega
}=-h_{0}$. Numerical results are similar to those for the Dirichlet-Neumann spectrum.

\section{Solution of the inverse problem}

Let us consider the inverse Sturm-Liouville problem consisting in recovering
the potential $q(x)$ from given several first eigenvalues of two spectra. For
definiteness we restrict ourselves to the Dirichlet-Dirichlet and
Dirichlet-Neumann spectra. Thus, given $\left\{  \mu_{n}^{2}\right\}
_{n=1}^{N_{1}}$ and $\left\{  \rho_{n}^{2}\right\}  _{n=0}^{N_{2}}$, first
eigenvalues of problems (\ref{SL equation}), (\ref{DD cond}) and
(\ref{SL equation}), (\ref{DN cond}), respectively. In the first step, as it
was explained in the previous section, we compute several coefficients
$\widetilde{s}_{n}(\pi)$, $n=1,\ldots,N$, $N\leq\min\left\{  N_{1}%
-1,N_{2}-1\right\}  $ and hence obtain the function $\widetilde{S}_{N}%
(\rho,\pi)$, which approximates the characteristic function of the problem
(\ref{Schrqtilde}), (\ref{DirProblem2}). In the second step we compute several
coefficients $\widehat{\sigma}_{0}(\pi)$,...,$\widehat{\sigma}_{N}(\pi)$ and
then complete the spectrum of problem (\ref{SL equation}), (\ref{DN cond}).
Next, let us consider the solution $\psi(\rho,x)$ of equation
(\ref{SL equation}) satisfying the initial conditions at $\pi$:
\[
\psi(\rho,\pi)=1,\quad\psi^{\prime}(\rho,\pi)=0.
\]
Analogously to the solution (\ref{phiNSBF}) the solution $\psi(\rho,x)$ admits
the series representation
\begin{equation}
\psi(\rho,x)=\cos\left(  \rho\left(  \pi-x\right)  \right)  +\sum
_{n=0}^{\infty}\left(  -1\right)  ^{n}\tau_{n}\left(  x\right)  \mathbf{j}%
_{2n}\left(  \rho\left(  \pi-x\right)  \right)  , \label{psi}%
\end{equation}
where $\tau_{n}\left(  x\right)  $ are corresponding coefficients, analogous
to $g_{n}\left(  x\right)  $ from (\ref{phiNSBF}). Similarly to
(\ref{qi from g0}) the equality
\begin{equation}
q(x)=\frac{\tau_{0}^{\prime\prime}(x)}{\tau_{0}(x)+1} \label{q from tau}%
\end{equation}
is valid.

Note that for $\rho=\rho_{k}$ the solutions $S(\rho_{k},x)$ and $\psi(\rho
_{k},x)$ are linearly dependent because both are eigenfunctions of problem
(\ref{SL equation}), (\ref{DN cond}). Hence there exist such real constants
$\beta_{k}\neq0$, that
\begin{equation}
S(\rho_{k},x)=\beta_{k}\psi(\rho_{k},x). \label{S=psi}%
\end{equation}
Moreover, these multiplier constants can be easily calculated by recalling
that $\psi(\rho_{k},\pi)=1$. Thus,
\[
\beta_{k}=S(\rho_{k},\pi)
\]
and we approximate these constants with the aid of the coefficients
$\widetilde{s}_{n}(\pi)$:%
\[
\beta_{k}\approx\widetilde{S}_{N}(\sqrt{\rho_{k}^{2}+\mu_{1}^{2}},\pi
)=\frac{\sin\left(  \sqrt{\rho_{k}^{2}+\mu_{1}^{2}}\pi\right)  }{\sqrt
{\rho_{k}^{2}+\mu_{1}^{2}}}+\frac{1}{\sqrt{\rho_{k}^{2}+\mu_{1}^{2}}}%
\sum_{n=0}^{N}(-1)^{n}\widetilde{s}_{n}(\pi)\mathbf{j}_{2n+1}(\sqrt{\rho
_{k}^{2}+\mu_{1}^{2}}\pi).
\]
Having computed these constants we use equation (\ref{S=psi}) for constructing
a system of linear algebraic equations for the coefficients $s_{n}(x)$ and
$\tau_{n}\left(  x\right)  $. Indeed, equation (\ref{S=psi}) can be written in
the form%
\begin{align}
&  \frac{1}{\rho_{k}}\sum_{n=0}^{\infty}(-1)^{n}s_{n}(x)\mathbf{j}_{2n+1}%
(\rho_{k}x)-\beta_{k}\sum_{n=0}^{\infty}\left(  -1\right)  ^{n}\tau_{n}\left(
x\right)  \mathbf{j}_{2n}\left(  \rho_{k}\left(  \pi-x\right)  \right)
\nonumber\\
&  =-\frac{\sin\left(  \rho_{k}x\right)  }{\rho_{k}}+\beta_{k}\cos\left(
\rho_{k}\left(  \pi-x\right)  \right)  . \label{system beta k}%
\end{align}

We have as many of such equations as many Dirichlet-Neumann singular numbers
$\rho_{k}$ are computed. For computational purposes we choose some natural
number $N_{c}$ - the number of the coefficients $s_{n}(x)$ and $\tau
_{n}\left(  x\right)  $ to be computed. More precisely, we choose a
sufficiently dense set of points $x_{m}\in(0,\pi)$ and at every $x_{m}$
consider the equations
\begin{align*}
&  \frac{1}{\rho_{k}}\sum_{n=0}^{N_{c}}(-1)^{n}s_{n}(x_{m})\mathbf{j}%
_{2n+1}(\rho_{k}x_{m})-\beta_{k}\sum_{n=0}^{N_{c}}\left(  -1\right)  ^{n}%
\tau_{n}\left(  x_{m}\right)  \mathbf{j}_{2n}\left(  \rho_{k}\left(  \pi
-x_{m}\right)  \right) \\
&  =-\frac{\sin\left(  \rho_{k}x_{m}\right)  }{\rho_{k}}+\beta_{k}\cos\left(
\rho_{k}\left(  \pi-x_{m}\right)  \right)  .
\end{align*}
Solving this system of equations we find $s_{0}(x_{m})$ and $\tau_{0}\left(
x_{m}\right)  $ and consequently $s_{0}(x)$ and $\tau_{0}\left(  x\right)  $
at a dense set of points of the interval $(0,\pi)$. Finally, with the aid of
(\ref{qi from s0}) or (\ref{q from tau}) we compute $q(x)$.

\section{Numerical examples}

Consider the potential from Example 2. In Figure 7 we show the recovered
potential in the case of ten pairs of the eigenvalues given. Here $N=9$, the
parameter $\omega$ was recovered with the accuracy $0.092$, and additionally
to the ten Dirichlet-Neumann eigenvalues given ninety eigenvalues were
computed by the spectrum completion technique. Their use when constructing the
system (\ref{system beta k}) is crucial, because the accuracy without these
additional eigenvalues computed is considerably worse.%

\begin{figure}
[ptb]
\begin{center}
\includegraphics[
height=3.0043in,
width=3.9967in
]%
{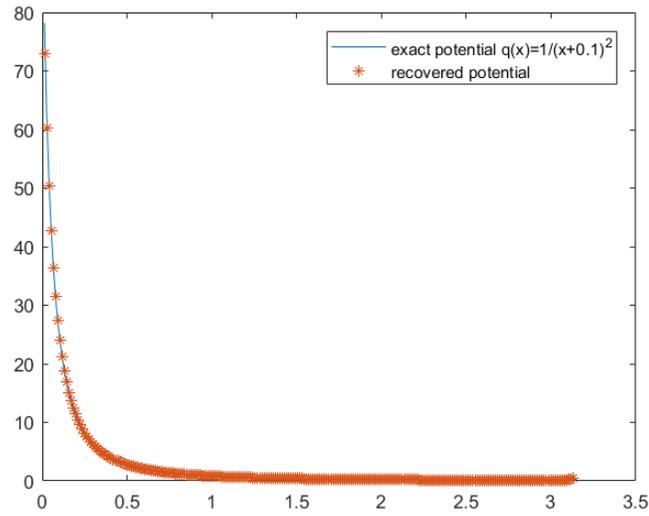}%
\caption{Potential from Example 2 recovered from 10 eigenpairs.}%
\label{Fig7}%
\end{center}
\end{figure}

\textbf{Example 3. }Consider a less smooth potential
\[
q(x)=\left\vert x-1\right\vert +1.
\]

In Fig. 8 we show the results of two spectra completion in the case of seven
eigenpairs given, $N=6$.%

\begin{figure}
[ptb]
\begin{center}
\includegraphics[
height=3.0636in,
width=4.0761in
]%
{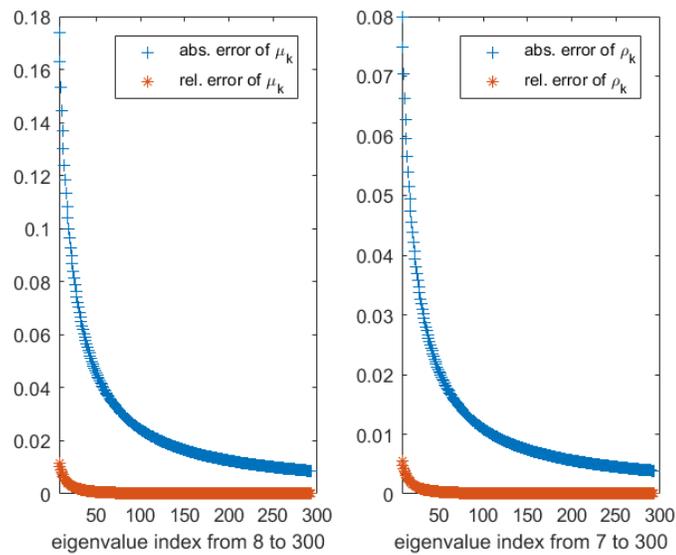}%
\caption{For the potential from Example 3, the accuracy of the
Dirichlet-Dirichlet (left) and Dirichlet-Neumann (right) eigenvalues computed
with the spectrum completion technique from seven eigenvalues of each spectrum
given.}%
\label{Fig8}%
\end{center}
\end{figure}

The result of the recovery of the potential in this case is presented in Fig. 9.%

\begin{figure}
[ptb]
\begin{center}
\includegraphics[
height=3.027in,
width=4.0273in
]%
{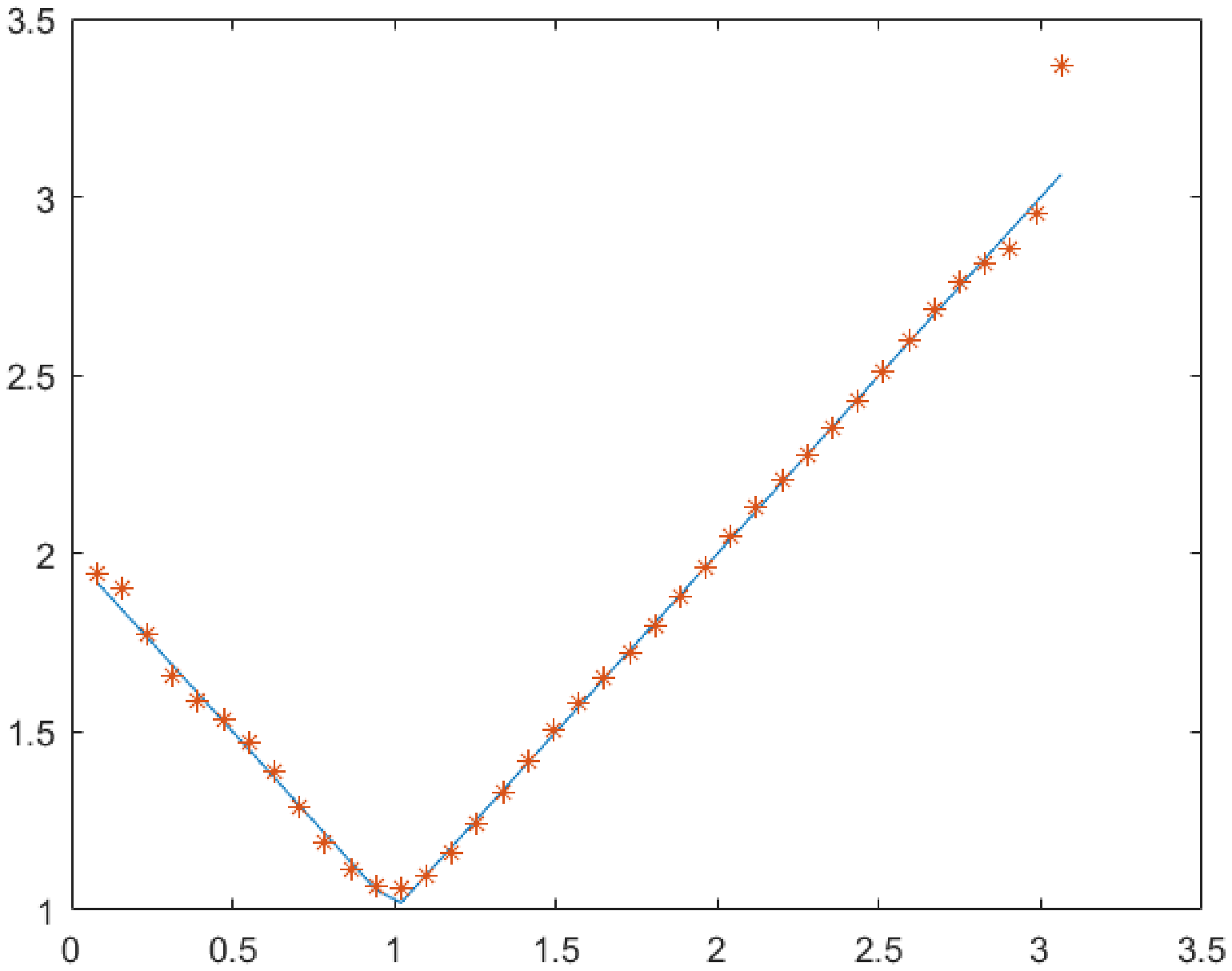}%
\caption{Potential from Example 3 recovered from seven eigenpairs given.}%
\label{Fig9}%
\end{center}
\end{figure}

Doubling the number of the eigenpairs given (for the same number of the
coefficients used, $N=6$) delivers more accurate results shown in Figs. 10 and 11.%

\begin{figure}
[ptb]
\begin{center}
\includegraphics[
height=3.2389in,
width=4.3089in
]%
{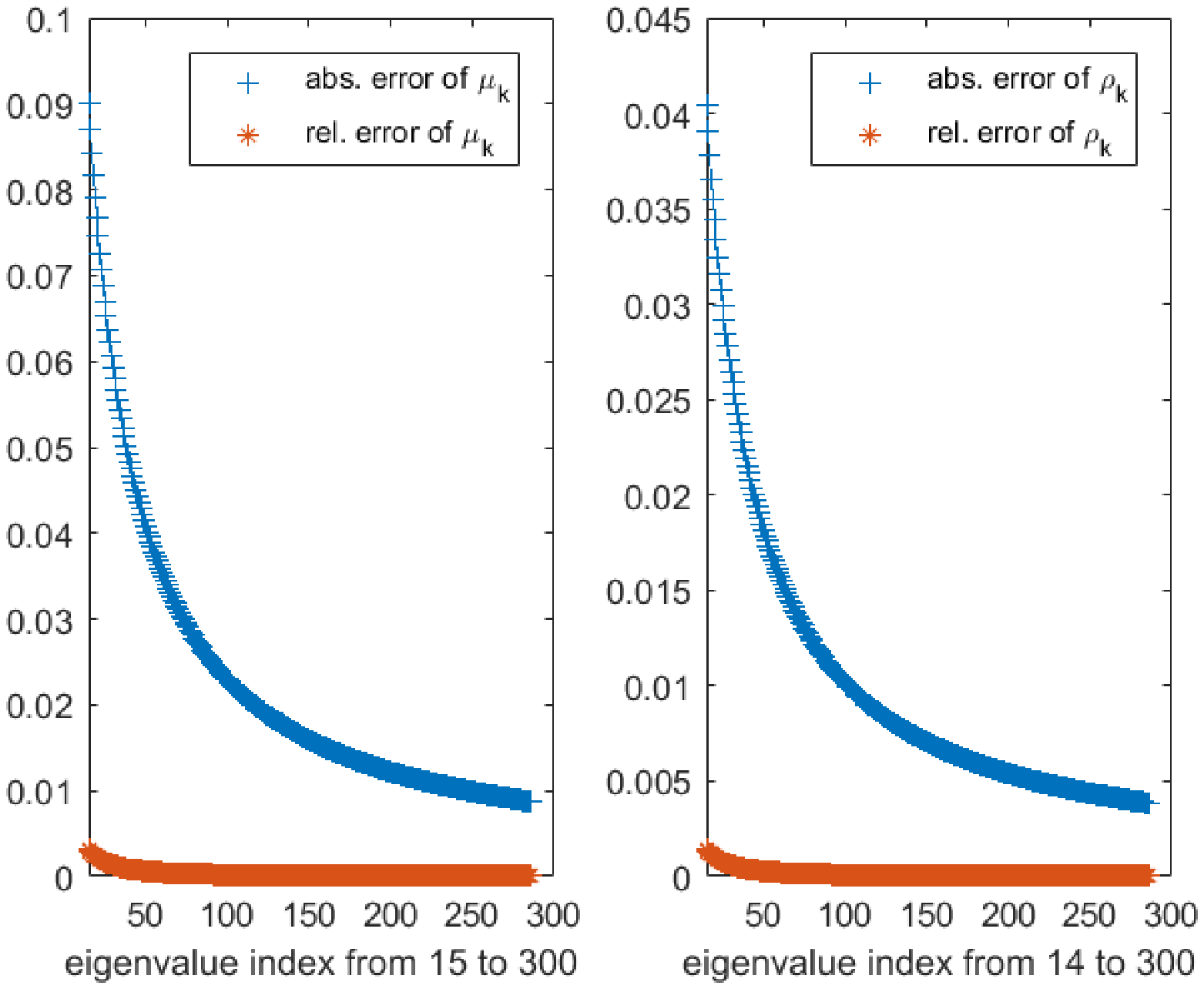}%
\caption{For the potential from Example 3, the accuracy of the
Dirichlet-Dirichlet (left) and Dirichlet-Neumann (right) eigenvalues computed
with the spectrum completion technique from 14 eigenvalues of each spectrum
given.}%
\label{Fig10}%
\end{center}
\end{figure}
%

\begin{figure}
[ptb]
\begin{center}
\includegraphics[
height=3.4281in,
width=4.561in
]%
{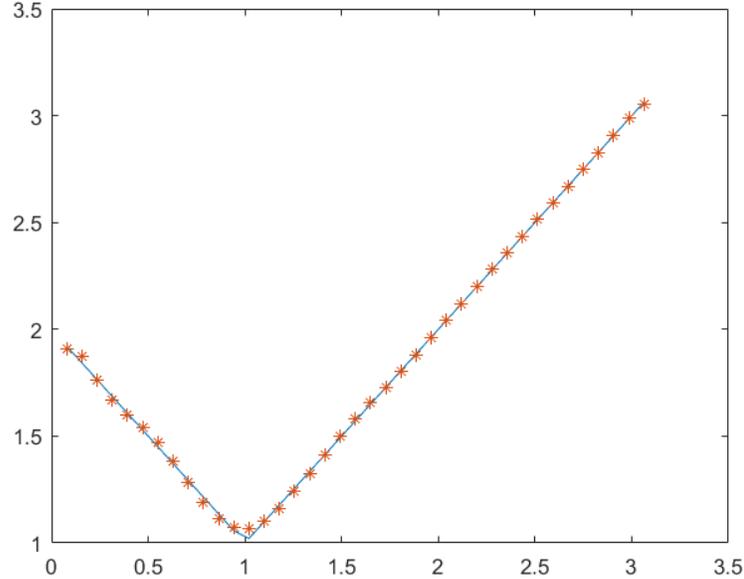}%
\caption{Potential from Example 3 recovered from 14 eigenpairs given. }%
\label{Fig11}%
\end{center}
\end{figure}

\section{Conclusions}

A simple method for completing the sequence of the eigenvalues of a regular
Sturm-Liouville problem is developed, which does not require neither
information on the potential nor the knowledge of the boundary conditions. The
spectrum is completed with a uniform absolute accuracy. Based on this spectrum
completion technique a direct method for solving two-spectra inverse
Sturm-Liouville problems on a finite interval is developed. The main role in
the proposed approach is played by the coefficients of the Neumann series of
Bessel functions expansion of solutions of the Sturm-Liouville equation and of
their derivatives. The given spectral data leads to an infinite system of
linear algebraic equations for the coefficients, and the potential is
recovered from the first coefficient alone.

The method is simple, direct and accurate. Its performance is illustrated by
numerical examples.

\section*{Acknowledgements}

Research was supported by CONACYT, Mexico via the project 284470 and partially
performed at the Regional mathematical center of the Southern Federal
University with the support of the Ministry of Science and Higher Education of
Russia, agreement 075-02-2022-893.

\end{document}